\documentclass{amsart}

\usepackage{array}
\usepackage{amssymb}
\usepackage{amsmath}
\usepackage{amsthm}
\usepackage{graphicx}
\usepackage{bm}
\usepackage{pgfplots}
\usepackage{caption}
\usepackage{algorithm}
\usepackage{algorithmic}
\usepackage{url}
\usepackage{longtable}
\usepackage{comment}
\usepackage[section]{placeins}

\newcolumntype{H}{>{\setbox0=\hbox\bgroup}c<{\egroup}@{}}

\newtheorem{theorem}{Theorem}

\setlongtables

\newcommand{\hl}[1]{{#1}}

\theoremstyle{definition}

\usepackage[square]{natbib}
\setcitestyle{square,semicolon,aysep={},notesep={, }}

\begin{document}

\title{Numerical and Statistical Analysis of Aliquot Sequences}

\author{K. Chum}
\address{University of Calgary \\
                2500 University Drive NW,
                Calgary, Alberta \\ Canada T2N 1N4}
\email{kchum@ucalgary.ca}

\author{R.~K. Guy}
\address{University of Calgary \\
                2500 University Drive NW,
                Calgary, Alberta \\ Canada T2N 1N4}
\email{rkg@cpsc.ucalgary.ca}

\author{M.~J.\ Jacobson, Jr.}
\address{University of Calgary \\
                2500 University Drive NW,
                Calgary, Alberta \\ Canada T2N 1N4}
\email{jacobs@cpsc.ucalgary.ca}
\thanks{The third author's research is supported by NSERC Discovery Grant 
RGPIN-2016-04545.}

\author{A.~S.\ Mosunov}
\address{University of Waterloo \\
                200 University Ave W,
               Waterloo, Ontario \\ Canada N2L 3G1}
\email{amosunov@uwaterloo.ca}


\begin{abstract}
We present a variety of numerical data related to the growth of terms in 
aliquot sequences, iterations of the function $s(n) = \sigma(n) - n$. First, we 
compute the geometric mean of the ratio $s_k(n)/s_{k-1}(n)$ of $k$th iterates 
for $n \leq 2^{37}$ and $k=1,\dots,10.$  Second, we extend the computation of 
numbers not in the range of $s(n)$ (called untouchable) by \citet{pollack:sum} 
to the bound of $2^{40}$ and use \hl{these data} to 
compute the geometric mean of the ratio of consecutive terms limited to terms 
in the range of $s(n).$  Third, we give an algorithm to compute $k$-untouchable 
numbers ($k-1$st iterates of $s(n)$ but not $k$th iterates) along with some 
numerical data. Finally, inspired by earlier work of 
\citet{devitt:thesis}, we estimate the growth rate of terms in aliquot 
sequences 
using a Markov chain model based on data extracted from thousands of sequences.
\end{abstract}

\maketitle

\section{Introduction}

An \emph{aliquot sequence} is the iteration of the function $s(n)=\sigma(n)-n$, 
where $\sigma(n)$ is the sum of the divisors of $n$. \citet{catalan}, 
corrected by \citet{dickson}, conjectured that all aliquot sequences 
terminate.   On the other hand, \citet{guy:drivers}
conjectured that, starting with an even value of $n$, many, perhaps almost all, 
such sequences diverge.  Aliquot sequences may terminate by reaching a prime 
$p$, since $s(p)=1$, by reaching a perfect number, for example 
$s(8128)=8128$, by reaching an amicable pair, for example $s(1184)=1210$ and 
$s(1210)=1184$, or by arriving in a longer cycle such as 14316, 19116, 31704, 
47616, 
83328, 177792, 295488, 629072, 589786, 294896, 358336, 418904, 366556, 274924, 
275444, 243760, 376736, 381028, 285778, 152990, 122410, 97946, 48976, 45946, 
22976, 22744, 19916, 17716, and then 14316 again.

\citet{devitt:thesis} used the average order of consecutive terms in a 
sequence, $s(n)/n$, taken over
even values of $n$, namely $5\pi^2/24\,-\,1=1.0562$ as evidence that
most aliquot sequences diverge, seemingly providing evidence in favor
of Guy and Selfridge.  However, it is the geometric mean, rather than
the arithmetic mean, that is relevant here.  \citet{bosma:constant} have 
calculated this to be no bigger than $\mu = 
0.967\ldots < 1$, computed to 13 decimal digits $\lambda = \log \mu = 
\hl{-0.0332594808010}\ldots$ 
by \cite{pomerance:constants}. 
Furthermore, \citet{pomerance:first} proved that 
the geometric mean of  $s_2(n)/s(n)$, taken over even $n > 2$, is also
equal to $\mu$, where $s_k(n)$ denotes the $k$-th iterate of
\hl{$s(n)$}. These results favor Catalan and Dickson.

However, there is a possible flaw in this argument, namely that the range of 
$s(n)$ varies.  \citet{erdos2},
\citet{pomerance:range}, and 
\citet{ppt:sn}
have obtained partial 
results about this.  For example,
$$s(192=2^6\cdot 3)=s(304=2^4\cdot 19)=s(344=2^3\cdot 43)=s(412=2^2\cdot 103)=316$$
so that 316 should perhaps be given correspondingly more weight than values of $n$ not 
in the range of $s(n)$. These latter seem difficult to calculate (what are the 
solutions to $s(x)=n$\,?) but experimental evidence \citep{pomerance:erdos} 
suggests that as many as one-third of the even numbers are not in the range of 
$s(n)$.

As there are numbers not in the range of $s(n)$, aliquot sequences will tend to 
be tributary at the numbers that \emph{are} in the range. Do these numbers, 
which should be counted by repetition, tend to be abundant (with $s(n) > n$) or 
deficient (with $s(n) < n$)?  
\citet{erdos2,erdos3} and others have investigated this aspect. \hl{By focusing 
on numbers $n$ that are highly abundant, \cite{erdos2} proved that there are 
infinitely many abundant numbers
not in the range of $s(n)$. Note that a similar question about deficient numbers remains open}.
Further, Erd\H{o}s improved on H.~W.~Lenstra's result that there are 
arbitrarily long increasing aliquot sequences \citep{erdos3}.

Another possible flaw in the argument of \citet{bosma:constant} and 
\citet{pomerance:first} is that the function $s(n)$ tends to 
preserve certain 
divisibility properties of $n$.  \citet{guy:drivers} 
explored the phenomenon of guides and drivers, particular divisors that tend to 
persist in consecutive sequence terms.  Most of the more persistent drivers 
cause terms to be abundant. \citet{pomerance:constants} has recently 
calculated the aliquot constants for terms with some divisibility restrictions, 
proving that $\log \mu$ is equal to \hl{$-0.3384354384114\ldots,$} 
\hl{$-0.2412950555350\ldots,$} and \hl{$0.1747760939329\ldots$} when restricted to terms 
that are even and square-free, congruent to $2 \bmod 4$, and divisible by 4, 
respectively.  The facts that guides and drivers are almost all not 
square-free, and that $4$ is itself a commonly-occurring guide, suggest that 
there may be a possibility that at least some aliquot sequences diverge.

In this work, we explored numerically two main lines of inquiry.
The first was whether, in order to accurately capture the behavior of
aliquot sequences, the geometric mean of $s(n)/n$ needs to take the
variability of the range of $s(n)$ into account as opposed to being
calculated over all even $n,$ for example, as many of these values
can never occur in a sequence.  This can be studied in a number of ways. One 
possibility is to study properties of $k$th iterates as in 
\citep{pomerance:first}, as the geometric mean of $s_k(n)/s_{k-1}(n)$ does 
indeed only take into account quantities in the range of $s(n).$  Another 
possibility is to enumerate explicitly \emph{touchable} numbers (those that are 
in the range of $s(n)$) and use only these when computing the geometric mean of 
$s(n)/n.$  Our second investigation was to account for the effects of guides 
and drivers by conducting numerical investigations measuring data 
occurring in actual sequences similar to \citet{devitt:thesis}; this 
method not only captures the variability of the range of $s(n)$, but also takes 
into account the influence of other factors such as guides and drivers 
\citep{guy:drivers} that the previous two approaches do not.

Our goal with this paper is to present novel numerical data on the growth of 
terms in aliquot sequences using all three of these approaches.  First, in 
Section~\ref{sec:kiterate} we describe our computation of the geometric mean 
$\mu_k(X=2^{37})$ of $s_k(n)/s_{k-1}(n)$ taken over all even $n \leq 2^{37}$ 
with 
$s_k(n) > 0$ for $1 \leq k \leq 10$.  In Section~\ref{sec:k-untouchable}, we 
describe the enumeration of untouchable numbers (those not in the range of 
$s(n)$), extending the table of \citet{pollack:sum} from 
a 
bound of $10^{10}$ to $2^{40}$. We use \hl{these data} to compute 
$\mu_1(2^{40})$ 
taken over even $n$ that are in the range of $s(n).$  We also present, in 
Section~\ref{sec:k-untouchable} an algorithm to compute $k$-untouchable 
numbers, those numbers in the range of $s_{k-1}(n)$ but not $s_k(n),$ and list 
all $k$-untouchable numbers up to $10^7$ for $k = 2, 3, 4, 5$ and
to $10^8$ for $k = 2, 3.$ Finally, in Section~\ref{sec:Devitt} we extend 
Devitt's work \citep{devitt:thesis} by 
calculating the geometric mean of $s(n)/n$ obtained from $8000$ 
randomly-selected aliquot sequences with initial terms of various sizes that 
were run until termination or until a term exceeded $2^{288}$.

Our data on $k$-th iterates show that 
\mbox{$\mu_k(2^{37}) > 1$} for $k \geq 6,$ suggesting that terms with many 
successive preimages are more likely to increase.  However, our data also show 
that $\mu_k(X)$ tends to decrease as the bound $X$ increases. In view of 
Pomerance's conditional result stated in \citep[Theorem 2.4]{pomerance:first}, 
it is natural to conjecture that $\mu_k(X)$ approaches $\mu$ for all $k$. On 
the other hand, the empirical estimate we 
obtain by extending Devitt's work is greater than one, suggesting that further 
extensions to the analytic results of 
\citet{bosma:constant} and \citet{pomerance:first,pomerance:constants} that 
also account 
for 
the effect of guides and drivers would be of great interest.

\section{Geometric Means of $k$th Iterates\label{sec:kiterate}}

\citet{bosma:constant} proved that the geometric mean of 
$s(n)/n$ 
taken over even $n$ is equal to \mbox{$\mu = 0.967\ldots < 1$}, and 
$$\log \mu = \hl{-0.0332594808010}\ldots$$ 
was computed to 13 decimal digits of accuracy by 
\citet{pomerance:constants}. Further, 
\citet{pomerance:first} proved that the geometric mean of 
$s_2(n)/s(n)$, taken over even $n > 2$, is also equal to $\mu$. Both of these 
results give strong probabilistic evidence that most aliquot sequences 
converge. 

A natural question is whether the same mean also holds for $k$th iterates,
where $k > 2$.  In some sense, considering these quantities for larger values
of $k$ might give a more accurate picture of the average behavior of aliquot
sequences, as the ratio of successive terms is being measured further along in a sequence as opposed to the first two or three terms.  Motivated by this
question, we performed some numerical computations of the 
quantity
$$
\lambda_k(X) = \frac{1}{\# B_k(X)}\sum\limits_{n \in 
B_k(X)}\log\frac{s_k(n)}{s_{k-1}(n)},
$$
where
$$
B_k(X) = \left\{n \in \mathbb N \colon \textrm{$n \leq X$, $n$ is even and 
$s_k(n) > 0$}\right\}.
$$

To evaluate $\lambda_k(X)$, we start by computing $\sigma(n)$ for all $n$ such 
that $1 \leq n \leq X$ and storing the resulting values into a lookup table 
$\Sigma$. This was done by using the algorithm of \citet{moews}. For each even 
$n$, we compute the value $s(n) = \sigma(n) - n$ 
and store it into one of the sets $\mathcal S_{1}$ or $\mathcal L_{1}$, 
depending on whether $s(n) \leq X$ or not. We compute 
$\lambda_1(X)$ directly from the values stored in $\mathcal S_{1}$ \hl{and $\mathcal L_{1}.$}

In order to compute $\lambda_k(X)$ for $k > 1$, we need to determine $s_k(n) = 
s(s_{k-1}(n))$ for all even $n \leq X$. To do so, for each $m$ in $\mathcal 
S_{k-1}$ we find $\sigma(m)$ in our lookup table $\Sigma$. In order to 
determine $\sigma(m)$ for all $m$ in $\mathcal L_{k-1}$, we factor $m$ 
directly. In our implementation, several values of $m$ in 
$\mathcal L_{k-1}$ are factored in parallel. Given all the $\sigma(m)$ values, we compute $s(m) = \sigma(m) 
- m$ and, in the case that $s(m) \neq 0$, store it in either $\mathcal 
S_k$ or $\mathcal L_k$, depending on whether $s(m) \leq X$ or not.  From 
\hl{these data},
we obtain $s(m)/m = s_k(n)/s_{k-1}(n)$ for some even number $n$ such that 
$s_{k-1}(n) = m$, and the resulting values will range over all even $n \leq X$ 
with $s_k(n) > 0$. Finally, we compute $\lambda_k(X)$ by iterating over \hl{$\mathcal S_k$ and $\mathcal L_k$}. At this point, 
we can discard the sets $\mathcal S_{k-1}$ and $\mathcal L_{k-1}$ and use the sets $\mathcal S_k$ and $\mathcal L_k$ in 
conjunction with the procedure described previously to compute 
$\lambda_{k+1}(X)$.

Since $\sigma(n)$ for $1 \leq n \leq 2^{37}$ fits into \texttt{unsigned long} 
(8 bytes), the size of the table $\Sigma$ reaches $1$ TB in size, while the 
tables $\mathcal S_1$ and $\mathcal L_1$ reach $512$ GB. None of those can fit 
into memory, and so the calculations were handled by partitioning our data into 
several files, storing the files into the hard disk, and loading them 
one-by-one into the memory.

All the computations described in this section, as well as in 
Sections~\ref{sec:untouchable} and \ref{sec:k-untouchable}, were carried on 
WestGrid's 
supercomputer Hungabee, located at the University of Alberta \citep{hungabee}, 
Canada. Hungabee is a 16 TB shared memory system with 2048 Intel Xeon cores, 
2.67GHz each. Each user of Hungabee may request at most 8 GB of memory per 
core. Also, Hungabee provides a high performance 53 TB storage space, which 
allows us to write to multiple disks in parallel. All the data, as well as the 
programs used for its computation, are available from the authors upon request.

Define $s_0(n) := n$. Using Hungabee, we evaluated $\lambda_k(X)$ for $k = 1, 
2, \ldots, 10$ and $X = 2^{15}, 2^{16}, \ldots, 2^{37}$ (see 
\mbox{Table~\ref{tab:lambda_data}}). The timings of our computations for 
\mbox{$X = 2^{37}$} are as follows. The computation of $\lambda_k(X)$ for $k = 
1, 2, \ldots, 10$ was mostly done sequentially, processing each of 1024 files 
individually, and using 128 processors in parallel solely for the purpose of 
factoring numbers that are larger than $2^{37}$. The computation of 1024 files 
storing $\sigma(n)$ for all $n \leq 2^{37}$ took 4h 42m 8s of real time, using 
256 processors (total CPU time 7w 1d 3h 46m 8s). With the table of $\sigma(n)$ 
stored in memory, the evaluation of $\lambda_1(X)$ took 3h 4m 39s of total CPU 
time. No parallel computations were needed in this case. The computation of a 
single value $\lambda_k(X)$ for $k \geq 2$ requires approximately 4 times more, 
as the computation of $s(s_{k-1}(n))$ is not as straightforward, and when 
$s_{k-1}(n) > 2^{37}$ it must be factored directly in order to compute 
$\sigma(s_{k-1}(n))$. The computation of $\lambda_k(X)$ for $k = 2, 3, \ldots, 
10$ took 10h 4m 20s, and the entire computation of $\lambda_1(X), \ldots, 
\lambda_{10}(X)$ for $X = 2^{37}$ took 3d 21h 43m 38s of real time using 128 
processors. It is difficult to estimate the total CPU time due to the fact that 
most of the program is sequential, and only factorization of large numbers is 
carried in parallel.
\begin{table}[htb]
	\begin{center}
		\begin{tabular}{| c || c | c | c | c | c | c | c | c |}
			\hline
			$X$ & $2^{15}$ & $2^{20}$ & $2^{25}$ & $2^{30}$ & $2^{35}$ & $2^{37}$\\
			\hline
			\hline
			$\lambda_1(X)$ & $-0.03336$ & $-0.03326$ &  $-0.03326$
			& $-0.03326$  & $-0.03326$ & $-0.03326$\\
			\hline
			$\lambda_2(X)$ & $-0.09523$ & $-0.05338$  &  $-0.04273$
			& $-0.03910$  & $-0.03749$ & $-0.03706$\\
			\hline
			$\lambda_3(X)$ & $-0.06434$ & $-0.01980$  & $-0.01399$
			& $-0.01531$ & $-0.01763$ & $-0.01849$\\
			\hline
			$\lambda_4(X)$ & $-0.07394$ & $-0.01333$  &  $-0.00568$
			& $-0.00752$ & $-0.01081$ & $-0.01205$\\
			\hline
			$\lambda_5(X)$ & $-0.07451$ & $-0.00830$  &  $0.00296$ & $0.00132$
			& $-0.00259$ & $-0.00411$\\
			\hline
			$\lambda_6(X)$ & $-0.08764$ & $-0.00536$  &  $0.01496$ & $0.00772$
			& $0.00322$ & $0.00145$\\
			\hline
			$\lambda_7(X)$ & $-0.09902$ & $-0.00053$  & $0.01655$ & $0.01478$
			& $0.00977$ & $0.00779$\\
			\hline
			$\lambda_8(X)$ & $-0.10820$ & $0.00033$  & $0.02212$ & $0.02059$
			& $0.01515$ & $0.01297$\\
			\hline
			$\lambda_9(X)$ & $-0.11139$ & $0.00316$  & $0.02822$ & $0.02669$
			& $0.02091$ & $0.01854$\\
			\hline
			$\lambda_{10}(X)$ & $-0.11341$ & $0.00708$  & $0.032823$ & $0.03193$
			& $0.02592$ & $0.02339$\\
			\hline
		\end{tabular}
	\end{center}
	\caption{Values of $\lambda_k(X)$.\label{tab:lambda_data}}
\end{table}

We report that the geometric means $\mu_k(X) = e^{\lambda_k(X)}$ exceed $1$ for 
$X = 
2^{37}$ and $k = 6,7,8,9,10$ when averaged over all even $n$ such that $s_k(n) 
> 0$. Moreover, as $k$ increases, the geometric means grow, too. \hl{As the 
function $\lambda_k$ pre-selects those $n$ that have not
yet reached a prime, it is not surprising that as $k$ increases,
so does $\lambda_k(X)$. What is more interesting is that}, as 
$k$ remains fixed, the geometric means decrease with the growth of $X$, 
possibly approaching the geometric mean of $s(n) / n$.

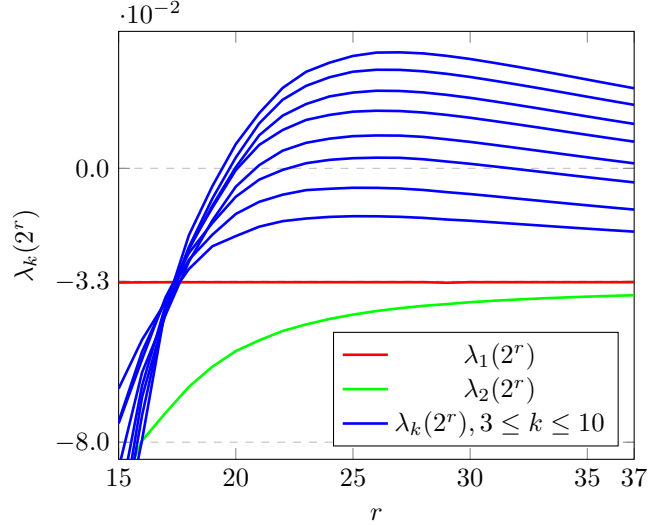
\begin{figure}[htb]
\centering
\captionof{figure}{Values of $\lambda_k(2^r)$ with $k$ fixed.\label{fig:x_varies}}
\begin{tikzpicture}
\begin{axis}[
    xlabel={$r$},
    ylabel={$\lambda_k(2^r)$},
    xmin=15, xmax=37,
    ymin=-0.085, ymax=0.040,
    xtick={15,20,25,30,35, 37},
    ytick={-0.08, -0.033, 0, 0.5, 0.1},
    legend pos=north west,
    ymajorgrids=true,
    grid style=dashed,
    legend pos=south east,
    y tick label style={
        /pgf/number format/.cd,
        fixed,
        fixed zerofill,
        precision=1,
        /tikz/.cd
    }
]


\addplot[
    line width=1.0pt,
    color=red,
    mark=none,
    ]
    coordinates {  
(15,-0.033360563730525767728488059729828975059)(16,-0.033310830683572155010794724332519567002)(17,-0.033274858288142921608963308179882256374)(18,-0.033272863364608875395546352772435086756)
(19,-0.033264928216387759835628873089393572336)(20,-0.033262938154635105381747427244152342845)(21,-0.03326080016534133099)(22,-0.03326034323393853992)(23,-0.03325989941606625339)(24,-0.03325970039985014920)(25,-0.03325955270086746601)(26,-0.03325953562779815353)(27,-0.03325950981718292510)(28,-0.03325949513604160857)(29,-0.03341203467096533902)(30,-0.03325948461190540023)(31,-0.03325948236600326235)(32,-0.03325948176194984396)(33,-0.03325948133096016945)(34,-0.03325948121093411092)(35,-0.03325948098624990218)(36,-0.03325948053271438748)(37,-0.03325948046763788774)
    };
    \addlegendentry{$\lambda_1(2^r)$}


\addplot[
    line width=1.0pt,
    color=green,
    mark=none,
    ]
    coordinates {  
(15,-0.095235281779911149461786136482365109075)(16,-0.079549094230760486015497213137612221799)(17,-0.071444948973371531667584809680792969117)(18,-0.063810591109932364158388699105590142295)(19,-0.057950380847866120786032237737301387272)(20,-0.053384053931977816201107135228753676969)(21,-0.05024607779847221206)(22,-0.04751094427445635288)(23,-0.04559803829877814157)(24,-0.04400793723155470294)(25,-0.04273437953813032097)(26,-0.04168900074390347738)(27,-0.04085568778591633365)(28,-0.04016684028730294065)(29,-0.03966650177282755180)(30,-0.03910395076572507206)(31,-0.03868502018643052431)(32,-0.03832514981659845138)(33,-0.03801072911999600623)(34,-0.03773311022348131705)(35,-0.03748642082339231213)(36,-0.03726512348544783887)(37,-0.03706496505472389041)
    };
    \addlegendentry{$\lambda_2(2^r)$}


\addplot[
    line width=1.0pt,
    color=blue,
    mark=none,
    ]
    coordinates {  
(15,-0.064339296473538842186612629992904719071)(16,-0.050072427927802953350642719018904340659)(17,-0.039032836340579313789853718943279779882)(18,-0.029400102661715211192171882980319486049)(19,-0.022776044930675331868068157485082166625)(20,-0.019798188523973043947835502458088876544)(21,-0.01704467177674217487)(22,-0.01543365783393091210)(23,-0.01466649037048324555)(24,-0.01419356249192458425)(25,-0.01398736711805494484)(26,-0.01401224165822368382)(27,-0.01417195495360873587)(28,-0.01447625675731249034)(29,-0.01495098106700092425)(30,-0.01530972762211209315)(31,-0.01577505708069932774)(32,-0.01624742697541741235)(33,-0.01671392833331293917)(34,-0.01718061119262993819)(35,-0.01763120581238645027)(36,-0.01807160925351631101)(37, -0.01849394208640171150)
    };
    \addlegendentry{$\lambda_k(2^r), 3 \leq k \leq 10$}


\addplot[
    line width=1.0pt,
    color=blue,
    mark=none,
    ]
    coordinates {  
(15,-0.073949295615029209279719979142034638409)(16,-0.053584731523885064304362428813763398836)(17,-0.038585969847145206902127342144642296256)(18,-0.026693700908853320339615144719403649082)(19,-0.019393779741977905883694821143274277348)(20,-0.013330064294402302835956530062359573388)(21,-0.00977844179678812216)(22,-0.00770153451685537133)(23,-0.00636733891370803474)(24,-0.00597412791044158245)(25,-0.00568187373706528342)(26,-0.00570871949592763328)(27,-0.00589232882804211454)(28,-0.00631796935770864883)(29,-0.00695744773187339377)(30,-0.00752319995015959721)(31,-0.00816459481992847305)(32,-0.00884379784734935434)(33,-0.00950432339300837960)(34,-0.01016728981637183768)(35,-0.01081240728751638007)(36,-0.01144045885755828364)(37,-0.01204620283054562598)
    };


\addplot[
    line width=1.0pt,
    color=blue,
    mark=none,
    ]
    coordinates {  
(15,-0.074506259471440746471849705102469811336)(16,-0.055766574098260926343203979442898769046)(17,-0.037722255886289652857048938678084820494)(18,-0.024727008235707322179612965994238143177)(19,-0.015966932812788626890091239542093908639)(20,-0.0083002592241293513713846665980651785540)(21,-0.00329595860018040785)(22,-0.00046056203975057758)(23,0.00149028001838813549)(24,0.00246995523275069353)(25,0.00295510060536073940)(26,0.00310949561701293889)(27,0.00306770133465552718)(28,0.00264991005749165044)(29,0.00198388126605934974)(30,0.00132247303977593928)(31,0.00056551779960030608)(32,-0.00022171862609593449)(33,-0.00101290376615679762)(34,-0.00181110269157039230)(35,-0.00259326877208793842)(36,-0.00336011796342954071)(37,-0.00410533547961502900)
    };


\addplot[
    line width=1.0pt,
    color=blue,
    mark=none,
    ]
    coordinates {  
(15,-0.087641649081984841112531608591060892012)(16,-0.059628520564613058004245371929528888776)(17,-0.042840673338217348547316309269824975841)(18,-0.027248406118546896725569305766222718268)(19,-0.015228312524061553848012138257828148298)(20,-0.0053608339237626686449914827972866389600)(21,0.00093621680158344352)(22,0.00515293119128809856)(23,0.00744832784094702147)(24,0.00866361810102312703)(25,0.00939226558992245150)(26,0.00959952209253823498)(27,0.00955170116426391248)(28,0.00917854515094941072)(29,0.00847234051537444459)(30,0.00772075834864260988)(31,0.00685333088479622363)(32,0.00595803754026093865)(33,0.00505014298495025382)(34,0.00412335550782348941)(35,0.00321626044955382285)(36,0.00232342541437224751)(37,0.00145360362585044605)
    };


\addplot[
    line width=1.0pt,
    color=blue,
    mark=none,
    ]
    coordinates {  
(15,-0.099020112628074058483861308611095104735)(16,-0.063497414037944643121731864500150104295)(17,-0.040071708343810143252232179548718849731)(18,-0.023217098783441599370219341071525186786)(19,-0.011505796061726486170133030066658501071)(20,-0.00052629532351183698275495280332776980617)(21,0.00684941769464097235)(22,0.01118534515547940644)(23,0.01413237734380834269)(24,0.01559455894718653375)(25,0.01654715868147051128)(26,0.01682155121429035338)(27,0.01672946801417224752
)(28,0.01630141663987220990)(29,0.01560236968299562148)(30,0.01477611757688332590)(31,0.01382060999276362931)(32,0.01284203536566951194)(33,0.01183072985344983327)(34,0.01079527469425888356)(35,0.00977486370967484662)(36,0.00877155388801996103)(37,0.00778928811245941169)
    };


\addplot[
    line width=1.0pt,
    color=blue,
    mark=none,
    ]
    coordinates {  
(15,-0.10820307950572049341534711642885189736)(16,-0.066812539848085028934955728984362545461)(17,-0.042420988268352768931166361054643963334)(18,-0.023626432904033262333223414433372651381)(19,-0.010933840275986950674878116209136337121)(20,0.00033017675058988258782052536281095975743)(21,0.00966476602539292039)(22,0.01538835217083454020)(23,0.01877784970846673243)(24,0.02099792176625295964)(25,0.02212468598777803108)(26,0.02266556333160705400)(27,0.02259157494473699665)(28,0.02217201840631182960)(29,0.02145894807115724875)(30,0.02058668750658789598)(31,0.01957158917835730150)(32,0.01850586008829469692)(33,0.01741109410322937145)(34,0.01627859033147934995)(35,0.01515343195346839471)(36,0.01404925336077354758)(37,0.01296524623530978700)
    };


\addplot[
    line width=1.0pt,
    color=blue,
    mark=none,
    ]
    coordinates {  
(15,-0.11139089662917225157394966109990069859)(16,-0.072159641411589498755605977995371756587)(17,-0.041192555725932857771156558581172796568)(18,-0.023002401365099149741843431329058701465)(19,-0.0088335325299543817974288308524292154137)(20,0.0031554304566617303026319808169877690668)(21,0.01312704751284974811)(22,0.02004256968606683495)(23,0.02408918061896641455)(24,0.02660172704264266580)(25,0.02821664011763224192)(26,0.02879670884919008703)(27,0.02875584711969692703)(28,0.02833466372099121763)(29,0.02760173416117642914)(30,0.02668622090168404695)(31,0.02560975684500050795)(32,0.02448963428689935815)(33,0.02332818536530424994)(34,0.02212990882960691463)(35,0.02091499713194492835)(36,0.01971781045933037810)(37,0.01854349636048996386)
    };


\addplot[
    line width=1.0pt,
    color=blue,
    mark=none,
    ]
    coordinates {  
(15,-0.11341166800534693542649846607698754205)(16,-0.079313583627745584376722032325919379227)(17,-0.040751866595721109201949564968148792701)(18,-0.019511331664494521129052336464341094391)(19,-0.0052116081510252686253097833732999154714)(20,0.0070756072016269640151885679747460199125)(21,0.01612977720668212420)(22,0.02348443500911108137)(23,0.02819549611690332774)(24,0.03086117020489615934)(25,0.03282304717818752132)(26,0.03380273323118635626)(27,0.03387642106889324906)(28,0.03356807489817675849)(29,0.03285360388158353739)(30,0.03193179896365901799)(31,0.03083304343622352006)(32,0.02967471509223966464)(33,0.02845770979869029216)(34,0.02720486428046264490)(35,0.02592183501009998406)(36,0.02464338929886494139)(37,0.02338735170118437512)
    };
\end{axis}
\end{tikzpicture}
\end{figure}

Figure~\ref{fig:x_varies} depicts $\lambda_k(2^r)$ as 
$r$ varies from $15$ to $37$ and $k$ remains fixed. The red line corresponds to 
$\lambda_1(2^r)$, the green line to $\lambda_2(2^r)$, and the blue lines, 
from bottom to top, correspond to $\lambda_3(2^r), \lambda_4(2^r), \ldots, 
\lambda_{10}(2^r)$. As the value of $X$ approaches infinity, the functions 
$\lambda_1(X)$ and $\lambda_2(X)$ converge to $\log \mu \approx -0.033$ from 
below. The convergence of $\lambda_1(X)$ is almost immediate and can hardly be 
observed. The convergence of $\lambda_2(X)$ is rather slow, and this was 
explained theoretically in \citep{pomerance:first}. As is seen in 
Figure~\ref{fig:x_varies}, the values of $\lambda_k(X)$ 
for $k = 3, 4, \ldots, 10$, unlike those for $k = 1, 2$, are strictly 
decreasing. We also see quite a remarkable difference between the behavior of
$\lambda_1(X)$, $\lambda_2(X)$ and $\lambda_k(X)$ for $k \geq 3$.

Let
$$
s^{-1}(n) := \{m \in \mathbb N : s(m) = n\}.
$$
It was conjectured by \hl{Erd\H{o}s, Granville, Pomerance, and Spiro see 
\citep[Conjecture 2.3]{pomerance:first}} that, if $A$ is a set of natural 
numbers of 
asymptotic density $0$, then $s^{-1}(A)$ also has asymptotic density $0$. 
Assuming this conjecture, Pomerance 
proved that for each integer $k \geq 2$ there is a subset $A_k \subseteq 
\mathbb N$  of asymptotic density $1$ with respect to $\mathbb N$ such that the 
geometric mean of $s_k(n)/s_{k-1}(n)$ averaged over all $n \in A_k$ converges 
to the Bosma-Kane constant $\mu$ \citep[Theorem 2.4]{pomerance:first}. Our 
computations serve as evidence of Pomerance's conditional result with $2\mathbb 
N$ in place of $\mathbb N$.

\hl{Finally, we have also evaluated $\lambda_{1}(X)$ up to $X=2^{40}$ in order 
to 
compare the value $-0.0332594808010\ldots$ computed to 
13 decimal digits of accuracy by \citet{pomerance:constants}. To 13 decimal 
places, we get $\lambda_1(2^{40}) = -0.0332594808048\ldots,$ 
matching Pomerance's computed value to the first 11 digits.}

\section{Tabulation of Untouchable Numbers\label{sec:untouchable}}

A number $n$ is called \emph{untouchable} if there is no $m$ such that $n = 
s(m)$. It is called \emph{touchable} otherwise. According to a slightly stronger variant of Goldbach's 
conjecture, for every odd integer $n \geq 9$ there exist distinct 
prime numbers $p$ and $q$ such that
$$
n = 1 + p + q = s(pq).
$$
The fact that $5$ is the only odd untouchable number follows from this 
conjecture, since $1 = s(2)$, $3 = s(4)$ and $7 = s(8)$, but no such expression 
exists for $n = 5$. This variant of Goldbach's conjecture has been verified 
computationally by Oliveira e Silva to $4 \cdot 10^{18}$ \citep{silva:goldbach}.

\citet{pollack:sum} gave a heuristic justification 
that the set of untouchable numbers has natural asymptotic density equal to
$$
\Delta = \lim\limits_{X \rightarrow \infty}\Delta(X),
$$
where
\begin{equation} \label{eq:untouchable_density}
\Delta(X) = \frac{1}{\log X}\sum\limits_{\substack{n \leq X\\2 \mid n}}\frac{1}{n}e^{-n/s(n)}.
\end{equation}
In order to verify the conjecture of Pollack and Pomerance computationally, we 
tabulated all even untouchable numbers to $2^{40}$. Previously, the tabulation 
of untouchable numbers was done by \citet{pomerance:erdos} 
to 
$10^8$ and by \citet{pollack:sum} to $10^{10}$. In all 
three cases the tabulation was done using the algorithm of Pomerance and Yang 
described in \citep{pomerance:erdos}. However, our computations were carried on 
a much larger scale, requiring some additional techniques, such as buffering 
and the usage of the hard disk. In the end, our counts of untouchable numbers
to $10^{10}$ matched those given by \citet{pollack:sum}.
The fourth author is grateful to Prof.\ Pollack for his help to verify the correctness
of our computations.

Our modification of the algorithm of Pomerance and Yang is as follows. In order 
to tabulate all untouchable numbers up to $X$, we determine all touchable 
numbers first, and then count the numbers that are left out. We use a bit array to represent each integer less than $X$, initialized to zero and split over $K$ files.  
We precompute $\sigma(n)$ for all odd 
$n$ in the range from $1$ to $X/2$ using the algorithm of \citet{moews}. Then 
we use the algorithm of Pomerance and Yang to tabulate 
all touchable numbers in parallel.  Each touchable number $m$ found is stored in one of $K$ buffers, each of which holds integers represented in a particular binary file.  When a buffer is full, the bit in the corresponding file that represents each integer stored in the buffer is set to one.  At the end, the zero bits in each file correspond to untouchable numbers.

Our computations were done in parallel, using 8 processors, with 1 thread per
processor. In order to maintain a reasonable balance between the total number
of files and the size of an individual file, the number of files was chosen to be equal
to 4096. With this setup, the program for tabulation of untouchable numbers to
$X = 2^{40}$ terminated in 3d 4h 7m 11s (the total CPU time is 3w 4d 8h 57m 28s). In order to reduce the number of disk I/O operations, as
well as to reduce the amount of storage space required, we compressed our data
by fitting the information on the type of every even number $n$ into 2 bits, depending
on whether the equation $n = s(m)$ has odd/even solutions, or both, or no solutions
at all. Thus, for example, the number $2$ has type \texttt{00}, since it is untouchable;
number $4 = s(9)$ has type \texttt{01}; number $10 = s(14)$ has type \texttt{10}; and
number $6 = s(6) = s(25)$ has type \texttt{11}. The size of our data is 128 GB.

\begin{table}[htb]
	\begin{center}
		\begin{tabular}{| r | r | r | r | r |}
		\hline
		$X$ & $U(X)$ & $U(X)/X$ & $\Delta(X)$ & $\Delta(X) - U(X)/X$\\
		\hline
		\hline
		$10^{4}$			& $1212$  & $0.121200$  & $0.161059$	& 0.039859\\
		\hline
		$10^{5}$			& $13863$ & $0.138630$  & $0.164577$	& 0.025947\\
		\hline
		$10^{6}$			& $150232$  & $0.150232$ & $0.166923$		& 0.016691\\
		\hline
		$10^{7}$			& $1574973$  & $0.157497$ & $0.168599$	& 0.011102\\
		\hline
		$10^{8}$			& $16246940$  & $0.162469$ & $0.169857$	& 0.007388\\
		\hline
		$10^{9}$		& $165826606$  &  $0.165827$ & $0.170834$	& 0.005007\\
		\hline
		$10^{10}$		& $1681871718$ & $0.168187$ & $0.171617$	& 0.003430\\
		\hline
		$10^{11}$		& $16988116409$ & $0.169881$ & $0.172257$	& 0.002376\\
		\hline
		$10^{12}$		& $171128671374$ & $0.171129$ & $0.172790$	& 0.001661\\
		\hline
		$2^{40}$		& $188206399403$ & $0.171173$ & $0.172810$	& 0.001637\\
		\hline
		\end{tabular}
	\end{center}
	\caption{Counts of untouchable numbers to $2^{40}$.\label{tab:untouchable}}
\end{table}

The counts are given in Table~\ref{tab:untouchable}. By $U(X)$ we denote the 
total count of untouchable numbers up to $X$. Our computations seem to suggest 
that the natural asymptotic density of the set of untouchable numbers exists. 
Looking back at \eqref{eq:untouchable_density}, for $X = 2^{40}$ the 
expression inside the limit has value $0.17281$, and this is quite close 
to our value $U(2^{40})/2^{40} \approx 0.17117$. The values of $\Delta(X)$ for 
different $X$ are also given in Table~\ref{tab:untouchable}. 
Figure~\ref{fig:untouchable} depicts the graphs of $U(X)/X$ and $\Delta(X)$. The
heuristics of Pollack and Pomerance suggests that the graphs of $U(X)/X$ and 
$\Delta(X)$ should
approach the same limiting value, and according to Figure~\ref{fig:untouchable} 
this certainly seems to be the case.

\input{plot_log_base2_untouchable.tex}

It was suggested by the second author, as well as by Carl Pomerance, to 
consider different variants of the quantity $\lambda_1(X)$ by replacing the set 
$B_1(X)$ with something else.  For example, averaging over the set of all even 
touchable numbers, rather than over all even numbers, would be more 
appropriate, as the untouchable numbers never really occur in aliquot sequences 
and therefore do not affect their behavior. According to 
Table~\ref{tab:untouchable}, over a third of all even numbers up to $2^{40}$ 
are 
untouchable, and so they \hl{might} influence the resulting quantity $\mu_1(X) = 
\exp(\lambda_1(X))$ quite significantly. We report that the geometric mean of 
$s(n)/n$ taken over touchable even $n \leq 2^{40}$ is equal to $0.915285$, which 
is less than $\mu \approx 0.967$. On the other hand, the geometric mean of 
$s(n)/n$ taken over untouchable even $n \leq 2^{40}$ is equal to $1.076$. Thus 
the computations suggest that the presence of untouchable numbers actually 
forces the geometric mean over all even numbers to increase rather than to 
decrease. 

\hl{Heuristically speaking, it is perhaps not so surprising that even numbers 
in the range tend to be somewhat less abundant than even numbers in the domain. 
Since any function $f(n)$ with $f(n)$ large compared with $n$ will miss many 
numbers, we do expect many $n$ not to appear in the range of $s(n)$ when $s(n)$ 
is significantly larger than $n$. In fact, this observation lies in the core of 
the argument of \citet{erdos2}.}

Another variant of $\lambda_1(X)$ that we considered incorporates the fact
that certain numbers have more preimages under $s(n)$ than others, and so
they could be influencing the behavior of aliquot sequences more. For example,
the number $160154$ has $15$ preimages under $s(n)$, namely
$$
\begin{array}{l l l l l l}
160154	& = s(152776)	& = s(183016)	& = s(251614)	& = s(260182)	& = s(296074)\\
		& = s(298234)	& = s(302842) 	& = s(310942)	& = s(311086)	& = s(312166)\\
		& = s(313822)	& = s(313966)	& = s(315082)	& = s(315226)	& = s(315262).
\end{array}
$$
Hence the value of $\log(s(n)/n)$ for $n=160154$ should perhaps be considered with the weight $15$
instead of weight $1$ as it was done in the computation of $\lambda_1(X)$. Clearly,
every untouchable number in this case would have weight $0$. A slight modification of our program for the tabulation of untouchable numbers allowed us to tabulate the values of $\#s^{-1}(n)$ for all even $n$ up to $X = 2^{40}$, and evaluate the quantity
$$
\tilde \lambda(X) = \left(\sum\limits_{\substack{n \leq X\\2 \mid n}}\#s^{-1}(n)\right)^{-1}\sum\limits_{\substack{n \leq X\\2 \mid n}}\#s^{-1}(n)\log\frac{s(n)}{n}.
$$

Once again, our computations were done in parallel, using 64 processors with 
one thread per processor. The number of files was chosen to be equal to 8192. 
With this setup, the program for the tabulation of $\#s^{-1}(n)$ for even $n$ 
up to $X = 2^{40}$ terminated in 1w 8h 34m 37s (the total CPU time is 1y 42w 3d 
20h 55m 28s). Evidently, in comparison to the tabulation of untouchable 
numbers, the tabulation of $\#s^{-1}(n)$ took significantly more time. The main 
reason for the degradation in performance lies in the fact that we recorded not 
the type of a number, but the actual count of the total number of preimages. 
While the type of a number fits into 2 bits, the information on the number of 
preimages fits into 8 bits,\footnote{We have $\#s^{-1}(690100611194) = 139$, 
and this is the maximum over all $\#s^{-1}(n)$ with even $n \leq 2^{40}$. For 
more details, see the OEIS sequence A283157 \citep{oeis}.} so the total number 
of disk I/O operations increased at least 4 times.
The total size of the data is 512 GB.

We report that for $X = 2^{40}$ the geometric mean $\exp(\tilde \lambda(X))$ is 
equal to $0.862346$, which is less than the Bosma-Kane constant $\mu$. In fact, 
it is even less than the geometric mean taken over touchable even numbers up to 
$X$.

\section{Tabulation of $k$-Untouchable Numbers\label{sec:k-untouchable}}

Let $k$ be a positive integer exceeding one. A \emph{$k$-untouchable number} is 
a number lying in the range of $s_{k-1}$, but not in the range of $s_k$. Just like
untouchable numbers, $k$-untouchable numbers do not appear in aliquot sequences
$\{s_i(n)\}_{i=0}^\infty$ for \mbox{$i \geq k$}, so they can also be 
disregarded when computing the geometric mean of $s(n)/n$. 

In this section, we
introduce an algorithm which tabulates all \mbox{$k$-untouchable} numbers up to $X$ and
$k \leq K$.
We use this algorithm to tabulate $k$-untouchable numbers to $10^7$ for $k = 2, 3, 4, 5$ and
to $10^8$ for $k = 2, 3$, and then use \hl{these data} to compute the geometric 
mean of
$s(n)/n$ over all even $n \leq X$ that are not $k$-untouchable for various values of $k$
and $X$.

Our computations were done as follows. First we precomputed the table of preimages for all even $n \leq X$ using the algorithm of Pomerance and Yang (here $X$ was chosen to be $2^{30}$). Then for each even $n$ we launched the recursive procedure described below (Algorithm~\ref{alg:R}) to search for preimages under $s_k(n)$. This procedure terminates if we find some even $m > X$ such that $n = s_i(m)$ with $i < k$. Our experimental observations suggest that when preimages for all even $n \leq X$ are precomputed, the algorithm produces all $k$-untouchable numbers up to $X/2^k$.

Though in practice we used this relatively straightforward approach, we also 
propose another algorithm that, unlike the one described previously, is suited 
for tabulation of all $k$-untouchable numbers up to a specified bound $X$. 
Algorithm~\ref{alg:k-untouchable} describes the tabulation 
procedure. In short, it works as follows: we maintain an array $T[X][K]$ 
with entries $T[n][k]$ equal to $1$ or $0$ depending on whether $n$ lies in 
the range of $s_k$ or not. With the help of the Pomerance-Yang algorithm, we 
initialize a reference table $P$ so that  $P[n] = \{m \colon n = s(m)\}$ for 
$n \leq X$; and a hash table $Q$ containing ordered triples 
$(n,m,1)$ such that $n \leq X$ and $n=s(m)$. We then proceed by calling the 
recursive subroutine $R(n,m,k)=R(n,m,k,K,X,P,Q,T)$ for each triple $(n,m,k)$ in 
$Q$ (for the implementation of $R(n,m,k)$ see Algorithm~\ref{alg:R}). The 
subroutine $R(n,m,k)$ iterates recursively over all 
possible preimages of $n$ under $s_{k'}$ for $k' \in \{k+1,\ldots, K\}$. It 
returns $0$ in one of the following three cases:
\begin{enumerate}
	\item if $R$ reaches $m$ such that $n = s_{k'}(m')$ for some odd $m'$ and 
	$k' \in \{k, k+1,\ldots, K\}$. In this case, by a slightly stronger version 
	of the Goldbach's conjecture described in the previous subsection, it must 
	be the case that $n$ is $k'$-touchable for each $k' \in \{k,k+1,\ldots,K\}$;
	\item if for some ordered triple $(n,m',k')$ our function $R$ reaches $k' = 
	K$, which means that $n = s_K(m')$, so $n$ is $k'$-touchable for each $k' 
	\in \{k,k+1,\ldots,K\}$;
	\item there exists $k' \in \{k,k+1,\ldots,K\}$ such that the identity 
	$n=s_{k'}(m')$ implies that $m'$ is untouchable. In other words, $n$ is 
	$(k'+1)$-untouchable.
\end{enumerate}

\begin{algorithm}[htb]
	\captionof{algorithm}{$R(n,m,k,K,X,P,Q,T)$ subroutine} 
	\begin{algorithmic}[1]
		\REQUIRE An ordered triple $(n,m,k)$ such that $n=s_k(m)$.
		\STATE $T[n][k] \leftarrow 1$
		\IF{$m$ is odd}
		\FOR{\textbf{all} $k' \in \{k+1, \ldots, K\}$}
		\STATE $T[n][k'] \leftarrow 1$
		\ENDFOR
		\ELSIF{$m > X$}
		\STATE $Q \leftarrow Q \cup \{(n, m, k)\}$
		\RETURN $m$
		\ELSIF{$k < K$}
		\STATE $M \leftarrow 0$
		\FOR{\textbf{each} $m' \in P[m]$}
		\STATE $V \leftarrow R(n,m',k+1,K,X,P,Q,T)$
		\IF{$V = 0$}
		\STATE \textbf{break}
		\ELSE
		\STATE $M \leftarrow \max\{M, V\}$
		\ENDIF
		\ENDFOR
		\RETURN M
		\ENDIF
		\RETURN 0
	\end{algorithmic}
	\label{alg:R}
\end{algorithm}

In all of the three cases, upon termination of $R(n,m,k)$ we remove all 
triples $(n',m',k')$ in $Q$ with $n' = n$, thus indicating that the number $n$ 
is completely processed and all the values of $T[n][k']$ for $k' \in 
\{1,2,\ldots, K\}$ are correct. If, however, the evaluation of $R(n,m,k)$ 
resulted in some positive integer $M$ instead of $0$, one must face the 
situation that each recursive call inside $R(n,m,k)$ eventually hits some 
triple $(n,m',k')$ with $k' \in \{k,k+1,\ldots,K-1\}$ such that $n=s_{k'}(m')$ 
and $m' > X$. While executing, the subroutine $R(n,m,k)$ populates $Q$ with 
such triples $(n,m',k')$ and upon termination returns the maximum $M$ over all 
occurring values of $m'$. After that, we call the Pomerance-Yang algorithm 
again and expand the reference table $P[n]$ up to $n \leq M$. We then call the 
function $R(n,m,k)$ for each $(n,m,k)$ in $Q$ again, but with $M$ in place of 
$X$. The call to $R(n,m,k)$ and the further expansion of $P$ are made as many 
times as needed until $Q$ is empty. At this point, the algorithm terminates.

\begin{algorithm}[htb]
	\captionof{algorithm}{Tabulation of $k$-untouchable numbers} 
	\begin{algorithmic}[1]
		\REQUIRE two positive integers $K > 1$ and $X > 1$.
		\ENSURE in an array $T[X][K]$, each entry $T[n][k]$ is equal to $0$ if 
		$n$ is $k$-untouchable, and to $1$ otherwise.
		\STATE
		\STATE $T[n][k] \leftarrow 0$ for all $2 \leq n \leq X$ and $1 \leq k 
		\leq K$
		\STATE $Q \leftarrow \varnothing$
		\STATE // Pomerance-Yang algorithm
		\FOR{$n \in \{2,3, \ldots, X\}$}
		\STATE $P[n] \leftarrow \varnothing$
		\FOR{\textbf{each} $m$ such that $n=s(m)$}
		\STATE $P[n] \leftarrow P[n] \cup \{m\}$
		\STATE $Q \leftarrow Q \cup \{(n,m,1)\}$
		\ENDFOR
		\ENDFOR
		\FOR{\textbf{each} $(n,m,k)$ in $Q$}
		\STATE $M \leftarrow R(n,m,k,K,X,P,Q,T)$
		\IF{$M = 0$}
		\STATE $Q \leftarrow Q \setminus \left\{(n', m', k') \colon 
		n'=n\right\}$
		\ELSE
		\STATE // Refine the reference table using the Pomerance-Yang algorithm
		\FOR{$n \in {X+2, \ldots, M}$}
		\FOR{\textbf{each} $m$ such that $n=s(m)$}
		\STATE $P[n] \leftarrow P[n] \cup \{m\}$
		\ENDFOR
		\ENDFOR
		\STATE $X \leftarrow M$
		\ENDIF
		\ENDFOR
	\end{algorithmic}
	\label{alg:k-untouchable}
\end{algorithm}

Let $U_k(X)$ denote the total number of $k$-untouchable numbers up to $X$. We computed $U_k(X)$ for $X=10^7$, \mbox{$k = 2, 3, 4, 5$}, and for $X = 10^8$, $k = 2, 3$. Note that our bound $10^8$ is significantly smaller than the bound $2^{40}$, which occurred in the tabulation of untouchable numbers and the values of $\#s^{-1}(n)$. The reason is that in this case we store not just the types of numbers and not just the total number of preimages, but the preimages themselves.

%
%
%
%

Our choice of parameters is as follows: up to $X = 2^{30}$, every even number 
$n$ has at most 64 preimages, each of which fit into \texttt{unsigned long} (8 
bytes). Thus the total size of our table of preimages consisting of 4096 files 
is $(2^{30} / 2) \cdot 64 \cdot 8$ bytes, which is equal to 256 GB. The 
tabulation of preimages was carried in parallel, using 64 processors, and the 
program terminated in 2h 20m 12s of real time (the total CPU time is 6d 5h 32m 
48s). Then we launched (a simple version of) the tabulation of $k$-untouchable 
numbers, which was done sequentially. The program for tabulation was launched 
twice. At the first launch, our program tabulated all $k$-untouchable numbers 
for $k = 1, 2, 3$ up to $134222590 > 10^8$, and it terminated in 7h 2m 32s. At 
the second launch, our program tabulated all $k$-untouchable numbers for $k = 
1, 2, 3, 4, 5$ up to $34100312 > 10^7$, and it terminated in 1m 23s.
\begin{table}[htb]
	\begin{center}
		\begin{tabular}{| r || r | r | r | r | r | r | r |}
		\hline
$X$		& $10^4$	& $10^5$		& $10^6$		& $10^7$		& $10^8$\\
		\hline
		\hline
$U_1(X)$ & $1212$	& $13863$	& $150232$	& $1574973$	& $16246940$\\
		\hline
$U_2(X)$	& $389$	& $4459$		& $50824$	& $554973$	& $5792792$\\
		\hline
$U_3(X)$	& $134$	& $1648$		& $19628$	& $217628$	& $2200638$\\
		\hline
$U_4(X)$	& $43$	& $594$		& $7110$		& $79387$	& \\
		\hline
$U_5(X)$	& $15$	& $208$		& $2408$		& $27913$	& \\
		\hline
		\end{tabular}
	\end{center}
	\caption{Counts of $k$-untouchable numbers to $10^8$.\label{tab:k-untouchable}}
\end{table}

Table~\ref{tab:k-untouchable} contains the counts of $k$-untouchable numbers. 
Table~\ref{tab:first_10_k-untouchable_numbers} contains first $10$ 
$k$-untouchable numbers for $k = 1, 2, 3, 4, 5$. The sequences of 
$k$-untouchable numbers for $k = 2, 3, 4, 5$ now appear in the Online 
Encyclopedia of Integer Sequences \citep{oeis} under identifiers A283152, 
A284147, A284156 and A284187, respectively.
\begin{table}[htb]
	\begin{center}
		\begin{tabular}{| r || r | r | r | r | r | r | r | r | r | r |}
		\hline
$k$		& $a_1$	& $a_2$	& $a_3$	& $a_4$	& $a_5$	& $a_6$	& $a_7$	& $a_8$	& $a_9$	& $a_{10}$\\
		\hline
		\hline
$1$		& $2$	& $5$	& $52$	& $88$	& $96$	& $120$	& $124$	& $146$	& $162$	& $188$\\
		\hline
$2$		& $208$	& $250$	& $362$	& $396$	& $412$	& $428$	& $438$	& $452$	& $478$	& $486$\\
		\hline
$3$		& $388$	& $606$	& $696$	& $790$	& $918$	& $1264$	& $1330$	& $1344$	& $1350$	& $1468$\\
		\hline
$4$		& $298$	& $1006$	& $1016$	& $1108$	& $1204$	& $1502$	& $1940$	& $2370$	& $2770$	& $3358$\\
		\hline
$5$		& $838$	& $904$	& $1970$	& $2176$	& $3134$	& $3562$	& $4226$	& $5038$	& $5580$	& $6612$\\
		\hline
		\end{tabular}
	\end{center}
	\caption{First $10$ $k$-untouchable numbers.\label{tab:first_10_k-untouchable_numbers}}
\end{table}

Define the quantity
\begin{equation} \label{eq:k-untouchable_average}
\hat \lambda_k(X) = \frac{1}{\#C_k(X)} \sum\limits_{n \in C_k(X)}\log\frac{s(n)}{n},
\end{equation}
where
$$
C_k(X) = \left\{n \in \mathbb N \colon \textnormal{\hl{$n \leq X$, $n$ is even and in the range of $s_k$}}\right\}.
$$
The value of $\hat \lambda_k(X)$ is equal to the average of $\log(s(n)/n)$ for 
all even $n \leq X$, with $1$-untouchable, $2$-untouchable, \ldots, 
$k$-untouchable $n$ disregarded. Table~\ref{tab:removing_k_untouchable_numbers} 
depicts the behavior of $\hat \lambda_k(X)$, and it is clear that removing more 
and more $k$-untouchable numbers seem to force the average to go down rather 
than up.
\begin{table}[htb]
	\begin{center}
		\begin{tabular}{| r || r | r | r | r | r | r |}
		\hline
		$X$		& $k = 0$		& $k = 1$		& $k = 2$		& $k = 3$		& $k = 4$	& $k = 5$\\
		\hline
		\hline
		$10^7$	& $-0.03326$	& $-0.07370$	& $-0.08077$	& $-0.08506$	& $-0.08618$	& $-0.08645$\\
		\hline
		$10^8$	& $-0.03326$	& $-0.07754$	& $-0.08982$	& $-0.09619$	&			&\\
		\hline
		\end{tabular}
	\end{center}
	\caption{Values of $\hat \lambda_k(X)$.\label{tab:removing_k_untouchable_numbers}}
\end{table}

\section{Discussion}
	
Although our results to this point seem to favor \hl{Catalan--Dickson}, they 
should 
still be taken with a grain of salt. Even the fact that the Bosma-Kane constant 
$\mu = 0.967\ldots$ is only slightly less than one suggests that there are 
many abundant numbers, and it could be the case that some aliquot sequences 
have a bias towards them. For example, let $S(X)$ and $NS(X)$ be the 
collections of all square-free and non-square-free even numbers up to $X$, 
respectively, and consider the quantities
$$
\lambda_S(X) = \frac{1}{\# S(X)}\sum\limits_{n \in S(X)}\log\frac{s(n)}{n},
$$
$$
\lambda_{NS}(X) = \frac{1}{\# NS(X)}\sum\limits_{n \in NS(X)}\log\frac{s(n)}{n}.
$$
We evaluated $\lambda_S(10^r)$ and $\lambda_{NS}(10^r)$ at $r = 3, 4, \ldots, 
10$, and the corresponding values are given in Table~\ref{tab:sqf_vs_nsqf}. Our 
computations show that $\lambda_S(X)$ approaches the limit 
$-0.3384354384093\ldots$ proved by \citet{pomerance:constants}, but 
also suggest that the limit
$$
\lim_{X \rightarrow \infty} \lambda_{NS}(X)
$$
exists and, as expected, is close to the limit $0.1747760953325\ldots$ for 
integers congruent to $0 \bmod{4}$ proved by 
\citet{pomerance:constants}.  Thus, most of square-free even numbers are 
deficient, while most of non-square-free even numbers are abundant. As it was 
remarked by the second 
author, the first fact is of no surprise, since the down-driver $2$ exactly 
divides every square-free even number. Also note that in this case the values 
of 
$\lambda_S(X)$ and $\lambda_{NS}(X)$ are not nearly as close to zero as the 
value of $\log \mu \approx -0.033$. Hence it could be that some aliquot 
sequences contain more non-square-free numbers than square-free numbers (such 
as 
those with the drivers of \citet{guy:drivers}), and 
because of this they are more likely to go to infinity. The famous aliquot 
sequences starting at $276$, $552$, $564$, $660$ and $966$, --- so-called 
Lehmer's Five, --- seem to support this heuristic. For example, out of the 
first 1651 terms of the aliquot sequence $\{s_k(276)\}$ only 596 terms are 
square-free; for $\{s_k(552)\}$, out of the first 982 terms only 28 are 
square-free; for $\{s_k(564)\}$, out of the first 3315 terms only 1157 are 
square-free; for $\{s_k(660)\}$, out of the first 827 terms only 65 are 
square-free; and for $\{s_k(966)\}$, out of the first 819 terms only 154 are 
square-free \citep{lehmersfive}.
\begin{table}[htb]
	\begin{center}
		\begin{tabular}{| r || r | r || r | r |}
		\hline
		$X$	& $\#S(X)$	& $\lambda_S(X)$	& $\#NS(X)$	& 
		$\lambda_{NS}(X)$\\
		\hline
		\hline
		$10^3$	& $204$	& $-0.3342131$	& $296$	& $0.1701133$\\
		\hline
		$10^4$	& $2027$	& $-0.3379401$	& $2973$	& $0.1740342$\\
		\hline
		$10^5$	& $20267$	& $-0.3384178$	& $29733$	& $0.1746938$\\
		\hline
		$10^6$	& $202640$	& $-0.3384442$	& $297360$	& $0.1747062$\\
		\hline
		$10^7$	& $2026416$	& $-0.3384393$	& $2973584$	& $0.1747115$\\
		\hline
		$10^8$	& $20264234$	& $-0.3384351$	& $29735766$	& $0.1747105$\\
		\hline
		$10^9$	& $202642377$	& $-0.3384354$	& $297357623$	& $0.1747109$\\
		\hline
		$10^{10}$	& $2026423710$	& $-0.3384354$	& $2973576290$	& $0.1747109$\\
		\hline
		\end{tabular}
	\end{center}
	\caption{Values of $\lambda_S(10^r)$ and $\lambda_{NS}(10^r)$ for $3 \leq r \leq 10$.\label{tab:sqf_vs_nsqf}}
\end{table}

\section{Empirical Estimate of Average Growth of Terms\label{sec:Devitt}}

One possible issue with the data and results described above is that they do not account for the behavior of terms in actual sequences. Thus, our final numerical experiments aim to estimate the geometric mean of $s(n)/n$ empirically.

Drivers and guides, in particular, are two aspects of aliquot sequences that 
the previous approaches do not take into account.
Sometimes a divisor of an integer $n$ may tend to persist throughout
repeated applications of the sum-of-proper-divisors function. If this
divisor happens to be abundant, then the sequence will tend to
increase so long as this divisor remains in subsequent terms. This phenomenon 
was captured in the notion of \emph{guides} and \emph{drivers}, developed by 
\citet{guy:drivers}.  A
\emph{guide} of $n$ is a divisor of $n$ consisting of the greatest power of 2 
that exactly divides $n$, i.e.\ $2^a||n$ along with any subset of the \hl{prime 
factors} of
$2^{a+1}-1$. If a guide is especially persistent,
it is called a \emph{driver}. Drivers are of the form $2^a v$, where $v$ is 
odd, $v\mid2^{a+1}-1$, and $2^{a-1} \mid \sigma(v).$
The following theorem enumerates all drivers.
\begin{theorem}[Theorem 2 of \citep{guy:drivers}]
	The only drivers are $2$, $2^3 \cdot 3$, $2^3 \cdot 3 \cdot 5$, $2^5 \cdot 
	3 \cdot 7$, $2^9 \cdot 3 \cdot 11 \cdot 31$, and the even perfect numbers.
\end{theorem}
Of these drivers, $2$ is a \emph{downdriver}, while $2 \cdot 3$, $2^2 \cdot 7$,
$2^4 \cdot 31$, and $2^5 \cdot 3 \cdot 7$ are \emph{updrivers}. Downdrivers 
cause an aliquot sequence to decrease, whereas updrivers cause sequences to
increase. We focus on these 24 guides: $2^0$, $2$, $2 \cdot 3$, $2^2$,
$2^2 \cdot 7$, $2^3$, $2^3 \cdot 3$, $2^3 \cdot 5$, $2^3 \cdot 3 \cdot 5$, 
$2^4$, $2^4 \cdot 31$, $2^5$, $2^5 \cdot 3$, $2^5 \cdot 7$, $2^5 \cdot 3 \cdot 
7$, $2^6$, $2^6 \cdot 127$, $2^7$, $2^7 \cdot 3$, $2^7 \cdot 5$, $2^7 \cdot 3 
\cdot 5$, $2^8$, $2^9$, and $2^{\geq 10}$. We say that a guide $g$ is \emph{in 
control of} a sequence of terms when $g$ is a guide of said terms.

\subsection{Description of Experiments}\label{sec:experiments}

In his 1976 M.Sc.\ thesis \citep{devitt:thesis}, Devitt presented 
theoretical and numerical evidence, using a ``new method of factoring called 
Pollard-Rho'', that the average order of $s(n)/n$ \hl{for $n$ even} is greater 
than 1.  Devitt's 
method was to construct a Markov chain with states corresponding to the 
smallest 24 guides.  Transition probabilities and averages of $s(n)/n$ for each 
state were estimated empirically by sampling terms from sequences with various 
sized terms.  

The advantage of this approach is that it captures data on actual sequences, including the effects of drivers and guides.  The disadvantage is that it is purely empirical.

We have repeated the experiments that Devitt performed in his thesis, but due 
to the greater computing power and faster integer factorization algorithms 
available today, to a much greater range.  In
addition, we used the geometric mean as opposed to the
arithmetic mean used by Devitt.

We collected data from 8000 aliquot sequences: eight sets of
1000 sequences, with each sequence starting at $2^{16+32n}+2k$, where
$0\leq n \leq 7$ and $0\leq k < 1000$.  Our goal in separating the
sequences into these eight \emph{stages} was to study termination
behavior as sequence terms get larger.  Each sequence was followed
until a term became greater than $2^{288}$.  \hl{We used the same upper bound 
for all stages, as this gives every sequence as much of a chance to acquire a 
down-driver (and perhaps terminate) as we could feasibly enable, and, in 
addition, more 
accurately accounts for the observed trend that most sequences tend to 
increase.}  We used 
the Aliqueit software 
\citep{aliqueit} to manage the factorization of sequence terms.  Aliqueit uses 
a 
combination of factorization packages to factor each term, 
including GMP-ECM, Msieve, Yafu, and GGNFS.  We used the ECPP \citep{ECPP} 
implementation of the elliptic curve primality proving algorithm 
\citep{Morain07} to provide 
rigorous primality proofs for all probable primes dividing any sequence term.

Data from every sequence term was recorded in matrices as follows.  \hl{The
rows and columns are both indexed by the 25 smallest guides listed in 
Table~\ref{tab:active-guides}, including $2^0=1$
as a first column --- a sump into which the terminating sequences go.  Let 
$g(n)$ denote the guide of the integer $n.$}
In the
first of a pair of matrices, for each $n$, is accumulated 1 in the
$g(s(n))$-th column and $g(n)$-th row.  In the second of each pair of
matrices is accumulated the \emph{amplification}, $\log s(n)/n$.

From each pair of matrices we deduced a probability matrix and an
amplification matrix.  Say that the first has accumulated $t_{i,j}$
terms and the second an amplitude of $a_{i,j}$.  Then the probability
matrix has entries $p_{i,j}=t_{i,j}/\sum_{i,j}t_{i,j}$.  The amplification 
matrix will
have entries $A_{i,j}=a_{i,j}/t_{i,j}$.  The \hl{observed average 
amplification} is
$$\mathcal{A}=\sum_{i,j}p_{i,j}A_{i,j};$$  
\hl{having $\mathcal{A}>0$, i.e., with $e^{\mathcal{A}}>1$ would favor 
Guy--Selfridge, whereas $\mathcal{A} <0$ would favor Catalan--Dickson}.
We recorded separate values of $\mathcal{A}$ for each stage, as well
as one for all sequences taken together.

In addition to studying the observed amplification values, we again followed 
Devitt and modeled
the sequences using Markov chains in an effort to deduce empirically
the expected behavior. This is done by repeatedly squaring the matrix
of transition probabilities (a stochastic matrix) until the entries
converge.  The resulting transition probabilities were used to compute
a value of $\mathcal{A}$ that captures more accurately what we would
expect on average.
	
No measures were taken to avoid the effects of tributaries on the
calculated values, but we expect the effects of these to be negligible for such large sequence terms.

\subsection{Data on Terminating Sequences}
	
In Table~\ref{tab:ends}, we give data on the number of observed occurrences of 
different types of termination.  From all $8000$ sequences, we observed that 
$1544$ reached a prime, $6392$ passed our limit, and $64$ entered a cycle, with 
a total of $5644436$ terms. No new cycles were discovered in the
process. Of those sequences that entered a cycle, only $7$ were perfect
numbers, $47$ were amicable numbers, $1$ reached a cycle with $5$ terms, and 
$9$ reached the cycle with $28$ terms.  Statistics corresponding to each of the 
eight stages are included in 
\hl{Table~\ref{tab:ends}}  The proportions of sequences which appear not to 
terminate increase
strikingly as the stage increases.
\begin{table}[htb]
\begin{center}
\begin{tabular}{r|cccccccc|c}
       & st 0 & st 1 & st 2 & st 3 & st 4 & st 5 & st 6 & st 7 & Overall \\ 
start  & $2^{16}$ & $2^{48}$ & $2^{80}$ & $2^{112}$ & $2^{144}$ & $2^{176}$ & 
$2^{208}$ & $2^{240}$ & \\
\hline
primes  &    709 &    282 &    183 &    122 &     87 &     76 &     43 &     42 
&  1544 \\
perfs   &      4 &      1 &        &        &        &      1 &      1 &       
&     7 \\
amics   &     26 &      7 &      3 &      6 &      3 &        &      1 &      1 
&    47 \\
cycles  &      6 &      2 &      1 &      1 &        &        &        &       
&    10 \\
\hline
\hline
ended   &    745 &    292 &    187 &    129 &     90 &     77 &     45 &     43 
&   1608 \\
open    &    255 &    708 &    813 &    871 &    910 &    923 &    955 &    957 
&   6392 \\
\% open & 25.5\% & 70.8\% & 81.3\% & 87.1\% & 91.0\% & 92.3\% & 95.5\% & 95.7\% 
& 79.9\% \\
\hline
\end{tabular}
\end{center}
\caption{Termination Statistics.\label{tab:ends}}
\end{table}

In Table~\ref{tab:ends_cong}, we give percentages of sequences that remained 
open (i.e. reached our upper bound of $2^{288}$) for sequences whose initial 
terms have various properties.  We give the percentages for all sequences, 
those whose initial terms are $0 \bmod 4$, those whose initial terms are $2 
\bmod 4$, those whose initial terms have the down-driver $2$ or the guide 
$2^3$, and those whose initial terms are either deficient or abundant.
\begin{table}[htb]
	\begin{center}
		\begin{tabular}{r|cccccccc|c}
			& st 0 & st 1 & st 2 & st 3 & st 4 & st 5 & st 6 & st 7 & Overall 
			\\ 
			start  & $2^{16}$ & $2^{48}$ & $2^{80}$ & $2^{112}$ & $2^{144}$ & 
			$2^{176}$ & 
			$2^{208}$ & $2^{240}$ & \\
			\hline
 			all         & 25.5\% & 70.8\% & 81.3\% & 87.1\% & 91.0\% & 92.3\% & 
 			95.5\% & 95.7\% & 79.9\% \\
			$0 \bmod 4$ & 26.6\% & 71.2\% & 83.0\% & 88.4\% & 91.8\% & 93.0\% & 
			96.2\% & 96.6\% & 80.9\% \\
			abundant    & 38.8\% & 76.2\% & 86.4\% & 90.7\% & 94.4\% & 
95.0\% & 97.3\% & 97.6\% & 84.6\% \\
			$2 \bmod 4$ & 24.4\% & 70.2\% & 79.6\% & 85.8\% & 90.2\% & 91.6\% & 
94.8\% & 94.8\% & 79.0\% \\
			deficient   & 12.6\% & 65.6\% & 76.4\% & 83.7\% & 87.6\% & 89.7\% & 
93.7\% & 93.8\% & 75.3\% \\
			$2$ or $2^3$&  9.5\% & 63.9\% & 73.4\% & 82.0\% & 87.2\% & 89.0\% & 
			93.0\% & 93.5\% & 73.9\% \\
			\hline
		\end{tabular}
	\end{center}
	\caption{Percentage of "Open" Sequences for Initial Terms with Different 
	Properties.\label{tab:ends_cong}}
\end{table}
We observe that these properties of the initial terms do have some effect on 
the eventual fate of a sequence.  Pomerance's results 
\citep{pomerance:constants} show that terms that are $0 \bmod 4$ are expected 
to 
increase on average whereas those that are $2 \bmod 4$ decrease, suggesting 
that sequences whose initial terms have these properties have a better chance 
to increase overall or terminate, respectively.  \hl{Similarly, sequences that 
start with a deficient number begin by decreasing, and hence may have a better 
chance to terminate, whereas those that begin with an abundant number should 
tend to increase more often.  Finally, sequences whose 
initial terms have as a guide the down-driver $2$ or $2^3$ should have an even 
better chance of terminating, as these terms should drive the sequence more 
persistently down.}  Our data shows that this is in fact what happens 
in practice.  However, the effect of these properties of the initial term 
diminishes as the starting term itself increases.

\subsection{Data on Non-Terminating Sequences}

In Table~\ref{tab:active-guides}, we list the guides that were active for $n$
when $s(n)$ breached $2^{288}$ for the $6392$ sequences that did not
terminate.

\begin{table}[htb]
\begin{center}
\begin{tabular}{r@{\hspace{20pt}}|r@{\hspace{20pt}}l}
	Guide & \# Active & Driver Type \\
\hline
                          $2^0$ &    0  & down-driver \\
                          $2^1$ &    0  & down-driver \\
                    $2 \cdot 3$ & 1124  & up-driver \\
                          $2^2$ &  374 \\
                  $2^2 \cdot 7$ & 1620  & up-driver  \\
                          $2^3$ &    8 \\
                  $2^3 \cdot 3$ &  584 \\
                  $2^3 \cdot 5$ &   75 \\
          $2^3 \cdot 3 \cdot 5$ &  773  & up-driver  \\
                          $2^4$ &  386 \\
                 $2^4 \cdot 31$ &  531  & up-driver  \\
                          $2^5$ &   12 \\
                  $2^5 \cdot 3$ &  158 \\
                  $2^5 \cdot 7$ &   52 \\
          $2^5 \cdot 3 \cdot 7$ &  278  & up-driver  \\
                          $2^6$ &  126 \\
                $2^6 \cdot 127$ &  143  & up-driver  \\
                          $2^7$ &    2 \\
                  $2^7 \cdot 3$ &   49 \\
                  $2^7 \cdot 5$ &   15 \\
          $2^7 \cdot 3 \cdot 5$ &   24 \\
                          $2^8$ &   33 \\
                          $2^9$ &   12 \\
$2^9 \cdot 3 \cdot 11 \cdot 31$ &    0 & up-driver  \\
                       $2^{10}$ &   13 \\
\hline
\end{tabular}
\end{center}
\caption{Active guides for non-terminating sequences.\label{tab:active-guides}}
\end{table}
None of the terms here having $2^9$ as a guide contained the driver
$2^9 \cdot 3 \cdot 11 \cdot 31$. It seems that a good portion of the terms are 
\hl{under the control of an up-driver}.  There are also surprisingly many terms 
that
have $2^2$ as a guide. Since Pomerance's result \citep{pomerance:constants} 
shows that the geometric mean of $\log(s(n)/n)$ for $n$ divisible by 4 is 
greater than zero, this is in fact further evidence in support of 
\hl{Guy--Selfridge}.

\subsection{Average Amplification}

For each guide that we considered, we recorded in Table~\ref{table:guidedata} 
the number of occurrences, the number of ``runs" 
(consecutive sequence terms of length $\geq 2$ with that guide), the average 
length of a run, the logarithm of the average amplification over all terms with 
that guide, and the average amplification over each run with that guide.  The 
initial term of each of the $8000$ sequences is excluded in these counts, as 
amplification is not defined for the initial terms that do not have a 
predecessor.  For the even guides, we also recorded the expected number of 
occurrences of each  guide, assuming that divisibility properties of the 
sequence terms are the same as for random integers.  For example, of the 
$5632116$ even sequence terms we would expect \hl{one third} of them, 
$1877372$, to contain a factor $2$ but no higher power of $2$. 
\begin{table}[htb]
\begin{center}
\begin{tabular}{ r  r  r  r  r  r }
Guide & Times & Expected & Runs & Average & Average \\
	  & Seen  & Number   &      & Length  & Amplification \\
\hline
$2^0$                 &    4320 &       - &  1777 &   2.431 & -0.777 \\
$2$                   & 1288576 & 1877372 & 42353 &  30.488 & -0.470 \\
$2 \cdot 3$           &  746474 &  938686 &  5080 & 147.207 &  0.226 \\
$2^2$                 &  858815 & 1206882 & 58467 &  14.718 &  0.069 \\
$2^2 \cdot 7$         &  937561 &  201147 &  9369 & 100.101 &  0.333 \\
$2^3$                 &  269723 &  375474 & 45753 &   5.907 & -0.032 \\
$2^3 \cdot 3$         &  213178 &  187737 & 12274 &  17.390 &  0.536 \\
$2^3 \cdot 5$         &  124637 &   93869 & 13619 &   9.162 &  0.322 \\
$2^3 \cdot 3 \cdot 5$ &  139573 &   46934 &  3177 &  43.953 &  0.802 \\
$2^4$                 &  316435 &  340652 & 65368 &   4.848 &  0.303 \\
$2^4 \cdot 31$        &  257444 &   11355 &  2051 & 125.529 &  0.362 \\
$2^5$                 &   83619 &  100574 & 32961 &   2.541 &  0.107 \\
$2^5 \cdot 3$         &   60418 &   50287 & 11673 &   5.182 &  0.640 \\
$2^5 \cdot 7$         &   40648 &   16762 &  4963 &   8.195 &  0.345 \\
$2^5 \cdot 3 \cdot 7$ &   54003 &    8381 &  1630 &  33.139 &  0.813 \\
$2^6$                 &   87858 &   87309 & 35140 &   2.504 &  0.357 \\
$2^6 \cdot 127$       &   57622 &     693 &   306 & 188.307 &  0.348 \\
$2^7$                 &   20953 &   23467 & 12059 &   1.740 &  0.096 \\
$2^7 \cdot 3$         &   14195 &   11734 &  5466 &   2.600 &  0.631 \\
$2^7 \cdot 5$         &    6941 &    5867 &  3410 &   2.038 &  0.431 \\
$2^7 \cdot 3 \cdot 5$ &    4511 &    2933 &  1407 &   3.212 &  0.909 \\
$2^8$                 &   23925 &   22000 & 13507 &   1.774 &  0.386 \\
$2^9$                 &   12880 &   11000 &  7831 &   1.647 &  0.399 \\
$2^{\geq 10}$         &   12127 &   11000 &  7852 &   1.545 &  0.391 \\
\hline
\end{tabular}
\end{center}
\caption{Data on occurrences, runs, and amplification for each 
guide.\label{table:guidedata}}
\end{table}

Although the down-driver $2$ was the most frequently occurring, the up-drivers 
collectively occurred much more frequently than the down-drivers.  Notice that 
the down-driver $2$ occurs \hl{only on about} $2/3$ of the number of occasions 
than 
would be expected, while the up-drivers $2^2\cdot7$, $2^3\cdot3\cdot5$, 
$2^4\cdot31$, $2^5\cdot3\cdot7$, $2^6\cdot127$, \ldots all occur several times 
more frequently than expected, suggesting that on average we might expect terms 
in aliquot sequences to increase rather than decrease.

Following the procedure described in Section~\ref{sec:experiments}, we computed 
the average amplification per term and the expected amplification via Markov 
chain analysis, individually for experiments at each of the eight stages and 
for all sequences collectively.  The results are given in 
Table~\ref{tab:avg-amp}.
\begin{table}[htb]
	\begin{center}
		\begin{tabular}{r|cccccccc|c}
			Amplification & st 0 & st 1 & st 2 & st 3 & st 4 & st 5 & st 6 & st 
			7 & Overall  
			\\ 
			\hline
			average  & $0.135$ & $0.127$ & $0.117$ & $0.111$ & $0.102$ & 
$0.086$ & $0.073$ & $0.048$ & $0.101$ \\
			expected & $0.127$ & $0.139$ & $0.136$ & $0.135$ & $0.132$ & 
$0.121$ & $0.116$ & $0.103$ & $0.125$ \\
			\hline
		\end{tabular}
	\end{center}
	\caption{Average Amplification $\log(s(n)/n)$ 
	Statistics.\label{tab:avg-amp}}
\end{table}
In all cases, both the observed average per term and the expected average are 
greater than zero, giving empirical evidence that terms in aliquot sequences 
grow on average.  \hl{The fact that the majority of sequence terms come from 
open sequences that are increasing certainly contributes to these statistics 
being greater than zero, but this is, based on our observations, an accurate 
assessment of the average behavior of a sequence.}
Overall, our results indicate that the geometric mean of 
$s(n)/n$ is approximately $e^{0.125} \approx 1.133,$ lending support to the 
\hl{Guy--Selfridge} conjecture.

\hl{In Table~\ref{tab:avg-amp-seeds}, we give the observed and expected 
averages when restricting to initial terms with particular properties.  These 
data are computed over all sequences considered from all 8 stages. 
\begin{table}[htb]
	\begin{center}
		\begin{tabular}{r|cccccc}
			Amplification & all & $0 \bmod 4$ & abundant & $2 \bmod 4$ & 
			deficient & $2$ or $2^3$ \\ 
			\hline
			average  & $0.101$ & $0.106$ & $0.122$ & $0.097$ & $0.084$ & 
			$0.080$ \\
			expected & $0.125$ & $0.127$ & $0.131$ & $0.124$ & $0.121$ & 
			$0.119$ \\
			\hline
		\end{tabular}
	\end{center}
	\caption{Average Amplification $\log(s(n)/n)$ 
		Statistics for Varying Initial Terms.\label{tab:avg-amp-seeds}}
\end{table}
As above, we expect that the initial term will have some influence over the 
average amplifications; initial terms that are $0 \bmod 4$ or abundant should 
cause a slight increase, whereas initial terms that are $2 \bmod 4,$ deficient, 
or starting with a down-driver should cause a decrease, with the largest 
occurring for the down-driver case.  This is exactly what we observed.  
However, it is noteworthy that the average amplification, in all cases, is 
still greater than zero, indicating that even when considering only sequences 
that start by decreasing, on average the terms will tend to increase, and 
lending 
even more support to the Guy--Selfridge conjecture.}

\section{Conclusions}

Our results are mixed, and do not lead to a definitive conclusion.  Most of 
the  approaches related to analytic methods suggest that terms in aliquot 
sequences tend to decrease on average, except when restricting to 
non-square-free even values.  On the other hand approximating the average 
growth 
of 
terms experimentally suggests that terms tend to increase.  We hypothesize that 
the discrepancy is due to the effect of guides and drivers, which are taken 
into 
account in the pure numerical estimates but not in those inspired by recent 
analytic results.  Extensions of the analytic methods that account for their 
effect, further extending the recent work of
\citet{pomerance:constants} would thus be a natural avenue for future work.

\section{Acknowledgements}
The authors wish to thank Carl Pomerance and the anonymous referee for many 
helpful comments and suggestions.

\bibliographystyle{apalike}

\end{document}